\documentclass[twocolumn, twoside, 10pt, a4paper]{IEEEtran}

\usepackage{graphicx,color,epsf}
\usepackage{amsmath,amsfonts,amssymb,amscd,bm}
\usepackage{tikz}
\usetikzlibrary{patterns}
\usepackage{pgfplots}
\pgfplotsset{compat=newest}
\usepackage{csquotes}
\usepackage{hyperref}
\usepackage{xcolor}
\usepackage{tudacolors}
\usepackage[english]{babel}
\usepackage{amstext}
\usepackage{mathtools}
\usepackage[doi=false,isbn=false,url=false,style=ieee,giveninits=true]{biblatex}
\usepackage{eso-pic}

\AddToShipoutPicture*{\footnotesize\sffamily\raisebox{1cm}{\hspace{1.65cm}\fbox{\parbox{\textwidth}{\copyright 2020 IEEE.  Personal use of this material is permitted.  Permission from IEEE must be obtained for all other uses, in any current or future media, including reprinting/republishing this material for advertising or promotional purposes, creating new collective works, for resale or redistribution to servers or lists, or reuse of any copyrighted component of this work in other works.}}}}

\begin{document}

\title{\Large \bf Parallel-in-Time Solution of Eddy Current Problems\\Using Implicit and Explicit Time-stepping Methods}

\author{I. Cortes Garcia$^{1,2}$ I. Kulchytska-Ruchka$^{1,2}$ M. Clemens$^{3}$ and S. Sch\"{o}ps$^{1,2}$
\medskip  \\ \affiliation
\thanks{Manuscript received xxx y, 20zz; revised xxx yy, 20zz and xxx 1, 20zz; accepted xxx 1, 20zz. Date of publication xxx yy, 20zz; date of current version xxx yy, 20zz. (Dates will be inserted by IEEE; published is the date the accepted preprint is posted on IEEE Xplore; current version is the date the typeset version is posted on Xplore). Corresponding author: F. A. Author (e-mail: f.author@nist.gov). Digital Object Identifier (inserted by IEEE).}}

\def\affiliation{
\normalsize
$^{1}$Technical University of Darmstadt, Computational Electromagnetics Group, Schlo\ss{}gartenstr. 8, 64289 Darmstadt, Germany\\
$^{2}$Technical University of Darmstadt, Centre of Computational Engineering,
Dolivostr. 15, 64293 Darmstadt, Germany\\
$^{3}$University of Wuppertal, Chair of Electromagnetic Theory, 
Rainer-Gr\"{u}nter-Str. 21, 42119 Wuppertal, Germany\\
E-mail: idoia.cortes@tu-darmstadt.de
\vspace{-5mm}
}

\markboth{CEFC Pisa 2020, 1c, 5g}{CEFC Pisa 2020, 1c, 5g}

\maketitle

\begin{abstract}
The time domain analysis of eddy current problems often requires the simulation of long time intervals, e.g. until a steady state is reached. 
Fast-switching excitations e.g. in pulsed-width modulated signals require in addition very small time step sizes that significantly increase
computation time. To speed up the simulation, parallel-in-time methods can be used. 
In this paper, we investigate the combination of explicit and implicit time integration methods 
in the context of the parallel-in-time method Parareal and using a simplified model for the coarse propagator. 
\end{abstract}

\section{Introduction}
The transient simulation of magnetoquasistatic fields on electric devices allows evaluating their behaviour and thus aids 
the design or optimisation process. 

\begin{figure}
	\centering
	\includegraphics[width = 0.3\textwidth]{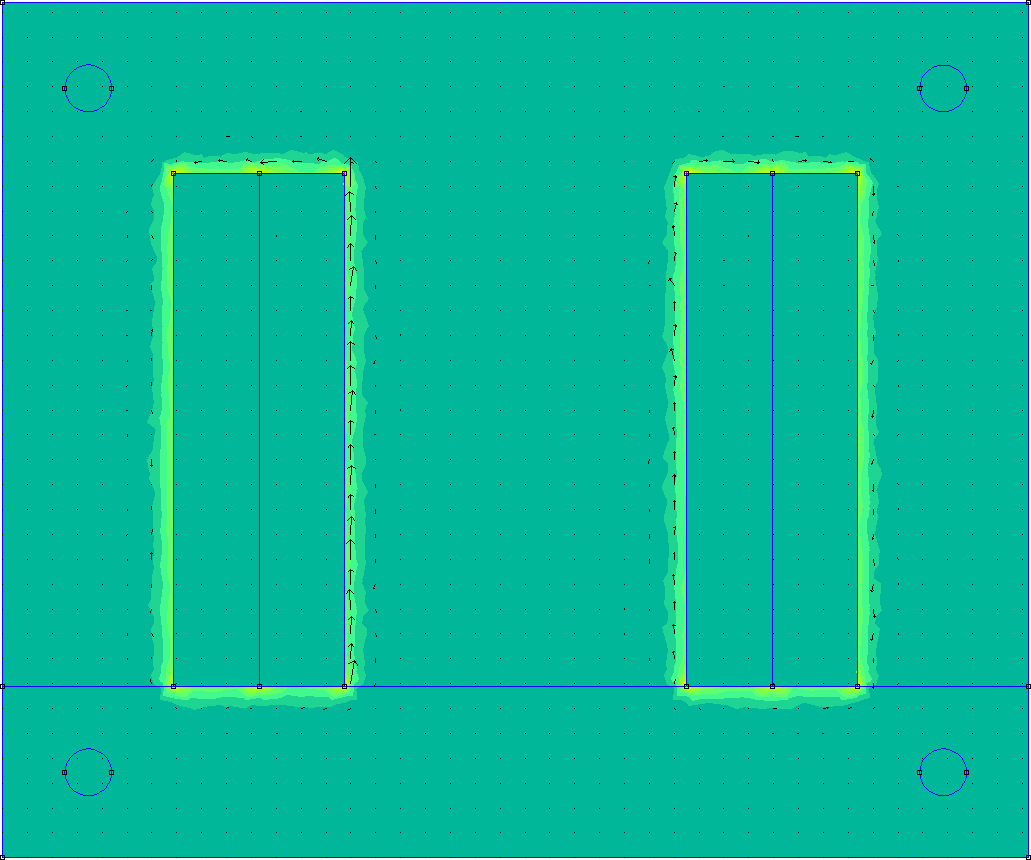}
	\caption{Transformer model of FEMM \cite{Meeker_2018aa}. }
	\label{fig:transformer}
\end{figure}

For the solution of the partial differential
equation, first a spatial discretisation is performed with e.g. the finite element method (FEM) \cite{Monk_2003aa} (see
Figure~\ref{fig:transformer}). This results in a
large system of only time dependent differential (algebraic) equations, that has to be integrated in time. Typically large system matrices and
 long time intervals yield computationally expensive problems, which become particularly inconvenient if fast-switching excitations are considered e.g.
 in a pulsed-width modulated (PWM) excitation scenario. Although the (magnetic) energy may vary slowly, the fast dynamics of the excitation require a fine temporal resolution and correspondingly small time steps. A reduction of this large computation time is possible by means of parallelisation.

 One approach towards shorter simulation times are
 domain decomposition methods. They allow the reduction of computation time by dividing the spatial domain, which in practice
 yields a set of smaller system matrices that can then be resolved in parallel. Parallel-in-time methods yield a second approach, 
 for example when
 domain decomposition methods are exploited up to their limit or 
 in cases where the time domain problem is the bottleneck of the simulation as is the case for example for the fast-switching excitations. 
 Recently, these
 type of methods have been proposed for the time-domain simulation of electric machines \cite{Takahashi_2019aa,Schops_2018aa}.
 
 Parareal is such a parallel-in-time method introduced in \cite{Lions_2001aa}, which can be interpreted as a special type of
 the multiple  shooting method
 \cite{Gander_2015aa}. In this article we propose a new combination of explicit and implicit time-integration methods in the context of 
 Parareal for discontinuous right-hand sides \cite{Gander_2019aa}.

The structure of the article is as follows:  Section~\ref{sec:eddy} presents the system of equations for the eddy current model
we consider. In Section~\ref{sec:time_int} the usage of explicit time integration methods for the spatially discretised eddy current equation
is discussed. Section~\ref{sec:Parareal} introduces
the Parareal algorithm and the newly proposed approach of combining it with implicit and explicit methods. Finally, Section~\ref{sec:numerics} presents
numerical simulations for a transformer model example and Section~\ref{sec:summary} closes with a summary and an outlook to future work.

\section{Eddy Current Model}\label{sec:eddy}
\begin{figure}
	\centering
	\includegraphics{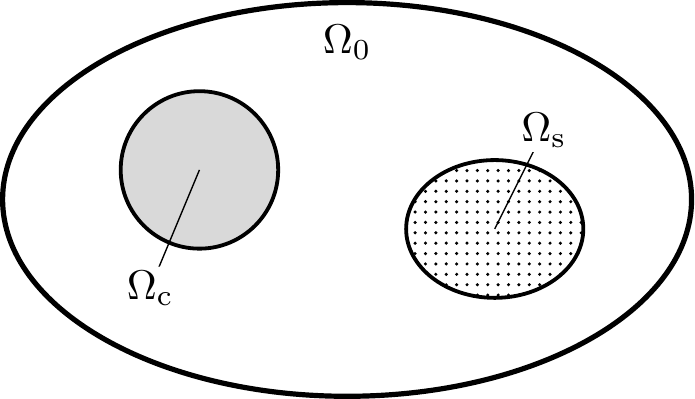}
	\caption{Sketch of the domain $\Omega$ where the eddy current partial differential equation is solved.}\label{fig:domain}
\end{figure}

We consider settings that can be described through a magnetoquasistatic approximation of Maxwell's equations  on domains
 $\Omega\subset\mathbb{R}^3$
and time intervals $\mathcal{I}=[t_0,\,t_{\mathrm{end}})\subset\mathbb{R}$.
The spatial domain $\Omega$ can be divided into three subdomains (see Figure~\ref{fig:domain}) 
$\Omega = \Omega_{\mathrm{s}}\cup\Omega_{\mathrm{c}}\cup\Omega_0$, where $\Omega_{\mathrm{s}}$ corresponds to the source domain,
 $\Omega_{\mathrm{c}}$ the domain of the conducting material and $\Omega_0$ the rest (typically the air region).
For these settings, the magnetic flux density $\vec{B}:\Omega\times\mathcal{I}\rightarrow\mathbb{R}^3$ can be described with the $\vec{A}^*$ 
formulation \cite{Emson_1988aa}. Here, a magnetic 
vector potential $\vec{A}:\Omega\times\mathcal{I}\rightarrow\mathbb{R}^3$ is defined, such that $\vec{B}=\nabla\times \vec{A}$ and the
electric field strength  is $\vec{E} = -\frac{\partial}{\partial t}\vec{A}$.
The resulting partial differential equation (PDE) that describes the magnetic flux density is the eddy current equation
\begin{equation*}
	\sigma\frac{\partial}{\partial t}\vec{A}+ \nabla\times(\nu\nabla\times \vec{A}) = \vec{\chi}_{\mathrm{s}}i(t)\;,
\end{equation*}
with $\vec{\chi}_{\mathrm{s}}:\Omega\times\mathcal{I}\rightarrow\mathbb{R}^3$ being a winding function that distributes
the source current on the spatial domain of the PDE \cite{Schops_2013aa} such that $\vec{\chi}_{\mathrm{s}}\mathbf{i}(t)$ is the
source current density.
The tensor
$\sigma:\Omega\rightarrow\mathbb{R}^{3\times 3}$ is the conductivity and $\nu:\mathbb{R}\times\Omega\rightarrow\mathbb{R}^{3\times 3}$
the possibly  nonlinear reluctivity. Note that the positive semidefinite
 conductivity $\sigma$ is only nonzero on the conducting region, that is,
$\mathrm{supp}\, \sigma = \Omega_{\mathrm{c}}$ and the positive definite  reluctivity $\nu$ is typically field-dependent on materials such as
iron, which are also part of $\Omega_{\mathrm{c}}$. The winding function $\vec{\chi}_{\mathrm{s}}$ is only prescribed at the source
domain $\Omega_{\mathrm{s}}$, that is, $\mathrm{supp}\, \vec{\chi}_{\mathrm{s}} = \Omega_{\mathrm{s}}$.

For the transient simulation of the eddy current problem typically the method of lines is used, where first a spatial
discretisation is performed and afterwards the time-dependent differential equations are solved with time integration
techniques. Let us consider the space-discrete eddy-current problem 
\begin{align}\label{eq:eddy}
\mathbf{M}_{\sigma}\frac{\mathrm{d}}{\mathrm{d}t}\mathbf{a}&= \underbrace{-\mathbf{K}_{\nu}(\mathbf{a})\mathbf{a}+
\mathbf{X}_{\mathrm{s}}\mathbf{i}(t)}_{\eqqcolon\mathbf{f}(\mathbf{a},t)}
\end{align}
with initial condition $\mathbf{a}(t_0)=\mathbf{a}_0$ on the interval $\mathcal{I}$ and conductivity matrix $\mathbf{M}_{\sigma}$, magnetic 
vector potential related degrees of freedom vector $\mathbf{a}(t)$, curl-curl-operator matrix $\mathbf{K}_{\nu}(\mathbf{a})$ and space-discrete 
winding function $\mathbf{X}_{\mathrm{s}}$. 
We assume the appropriate boundary conditions are contained in the system matrices. To ensure 
uniqueness of solution of \eqref{eq:eddy} in a three dimensional problem, an additional gauging condition has to be imposed. Here, for example, a tree-cotree gauging \cite{Munteanu_2002aa} or grad-div-regularisation \cite{Clemens_2002aa} can be employed. This ensures that the matrix pencil
$\lambda\mathbf{M}_{\sigma} + \mathbf{K}_{\nu}(\mathbf{a})$ is regular for $\lambda\in \mathbb{R}$ and thus the 
unique solvability of the (semi-discrete) system is guaranteed.

\section{Time Integration}\label{sec:time_int}
Time integration methods can be classified into explicit and implicit ones. The appropriate type of method for a specific system 
depends on properties such as its stiffness or nonlinearities. We consider the eddy current problem  \eqref{eq:eddy} and apply
two different time integration schemes to analyse their advantages and disadvantages.

As a first approach we consider the backward differentiation implicit Euler method. Its application to \eqref{eq:eddy} with a time step
size of $H$ yields at time instant $t_{i+1}$
\begin{equation*}
\mathbf{M}_{\sigma}\frac{\mathbf{a}_{i+1}-\mathbf{a}_{i}}{H} =-\mathbf{K}_{\nu}(\mathbf{a}_{i+1})\mathbf{a}_{i+1}+
\mathbf{X}_{\mathrm{s}}\mathbf{i}(t_{i+1})\;.
\end{equation*}
Whereas this method is unconditionally stable, which implies that the time step size $H$ has only to be reduced for accuracy reasons,
the system has to be resolved for $\mathbf{a}_{i+1}$. This involves, for nonlinear systems, the usage of a root-finding algorithm such as e.g. the
Newton method and repeated solutions of linear equation systems e.g. by factorization or iterative solvers. This increases the computational cost, as each time integration step requires several internal iterations of the
root-finding algorithm.

Due to the structure of \eqref{eq:eddy}, the usage of an explicit time integration method requires special treatment of the system of equations.
As the conductivity is zero in the non-conducting region $\Omega_{\mathrm{c}}^{\complement}$, the mass matrix $\mathbf{M}_{\sigma}$ is singular,
that is, the system of equations \eqref{eq:eddy} is a system of differential-algebraic equations (DAEs). If we follow the approach in \cite{Dutine_2017aa} and divide the degrees of freedom $\mathbf{a}$ into the ones corresponding to basis functions lying in the conducting region $\mathbf{a}_{\mathrm{c}}$ and the rest $\mathbf{a}_{\mathrm{nc}}$, then we can write system \eqref{eq:eddy} as
\begin{align*}
	\begin{pmatrix}
	\bar{\mathbf{M}}_{\sigma} & 0 \\
	0 & 0
	\end{pmatrix}
	\frac{\mathrm{d}}{\mathrm{d}t}
	\begin{pmatrix}
	\mathbf{a}_{\mathrm{c}}  \\
	\mathbf{a}_{\mathrm{nc}}
	\end{pmatrix} ={}&-
	\begin{pmatrix}
	\mathbf{K}_{\nu 1,1}(\mathbf{a}) & \mathbf{K}_{\nu 1,2} \\
	\mathbf{K}_{\nu 1,2}^{\top} & \mathbf{K}_{\nu 2,2}
	\end{pmatrix}
	\begin{pmatrix}
	\mathbf{a}_{\mathrm{c}}  \\
	\mathbf{a}_{\mathrm{nc}}
	\end{pmatrix}\\
	&+
	\begin{pmatrix}
	\bar{\mathbf{X}}_{\mathrm{s}}  \\
	0
	\end{pmatrix}\mathbf{i}(t)\;,
\end{align*}
where $\bar{\mathbf{M}}_{\sigma}$ is a positive definite matrix.  These type of systems of differential algebraic equations are often percieved as infinitely stiff problems \cite{Brenan_1995aa} and thus no explicit time integration methods can be employed. However, in \cite{Schops_2012aa, Dutine_2017aa} an approach is presented to circumvent this by means of a Schur complement. For a gauged system, due to the regularity of the matrix pencil
$\lambda\mathbf{M}_{\sigma} + \mathbf{K}_{\nu}(\mathbf{a})$, the matrix  $\mathbf{K}_{\nu 2,2}$ is regular. 
In this case, the Schur complement can be applied to extract an ordinary differential equation (ODE) for $\mathbf{a}_{\mathrm{c}}$
\begin{align}
\frac{\mathrm{d}}{\mathrm{d}t}	\mathbf{a}_{\mathrm{c}} ={}&- 
\bar{\mathbf{M}}_{\sigma}^{-1}(\mathbf{K}_{\nu 1,1}(\mathbf{a}_{\mathrm{c}}) - \mathbf{K}_{\nu 1,2}\mathbf{K}_{\nu 2,2}^{-1}\mathbf{K}_{\nu 1,2}^{\top})
\mathbf{a}_{\mathrm{c}}\nonumber \\
&-
\bar{\mathbf{M}}_{\sigma}^{-1}\mathbf{K}_{\nu 1,2}\mathbf{K}_{\nu 2,2}^{-1}\bar{\mathbf{X}}_{\mathrm{s}}\mathbf{i}(t)\;.\label{eq:ode}
\end{align}
Now the explicit Euler method can be employed on the derived ODE which yields for time $t_{i+1}$ and step size $h$ the expression
\begin{align*}
\frac{\mathbf{a}_{\mathrm{c},i+1}-\mathbf{a}_{\mathrm{c},i}}{h} ={}&- 
\bar{\mathbf{M}}_{\sigma}^{-1}(\mathbf{K}_{\nu 1,1}(\mathbf{a}_{\mathrm{c},i}) - \mathbf{K}_{\nu 1,2}\mathbf{K}_{\nu 2,2}^{-1}\mathbf{K}_{\nu 1,2}^{\top})
\mathbf{a}_{\mathrm{c},i}\nonumber \\
&-
\bar{\mathbf{M}}_{\sigma}^{-1}\mathbf{K}_{\nu 1,2}\mathbf{K}_{\nu 2,2}^{-1}\bar{\mathbf{X}}_{\mathrm{s}}\mathbf{i}(t_i)\;.
\end{align*}
To obtain the magnetic vector potential on the nonconducting region, the following equation can be exploited
\begin{equation*}
\mathbf{a}_{\mathrm{nc},i+1} = - \mathbf{K}_{\nu 2,2}^{-1}\mathbf{K}_{\nu 1,2}^{\top}\mathbf{a}_{\mathrm{c},i+1}\;.
\end{equation*}
Here, the time step size influences the accuracy of the solution, but also the stability of the time integration method. Thus, a sufficiently
small step size $h$ has to be chosen, to ensure the integration scheme remains within its stability region \cite{Hairer_1996aa}. This can
significantly reduce the required step size $h$ and therefore increase computational cost. However,
in contrast to the implicit methods, no (non)linear systems have to be resolved to obtain the solution for $\mathbf{a}_{\mathrm{c},i+1}$.

Note that the system contains two inverse matrices $\bar{\mathbf{M}}_{\sigma}^{-1}$ and $\mathbf{K}_{\nu 2,2}^{-1}$, that are in practice
not computed explicitly, but the corresponding linear system is solved.  This two operations can 
be computed efficiently as both matrices are constant and do not depend on the solution. Therefore, efficient techniques such as mass lumping on $\bar{\mathbf{M}}_{\sigma}$ and e.g. an LU decomposition of the linear matrix  $\mathbf{K}_{\nu 2,2}$ can be employed.
Furthermore, the expression 
$$\bar{\mathbf{Y}_{\mathrm{s}}}\coloneqq\bar{\mathbf{M}}_{\sigma}^{-1}\mathbf{K}_{\nu 1,2}\mathbf{K}_{\nu 2,2}^{-1}\bar{\mathbf{X}}_{\mathrm{s}}$$
on \eqref{eq:ode} has only to be computed once at the beginning of the simulation and can be interpreted as a different type of (constant) 
winding function
that is multiplied by the time dependent current $\mathbf{i}(t)$.

\section{Parareal}\label{sec:Parareal}
Parareal is a parallel-in-time method that takes advantage of parallel hardware by splitting the time interval
$\mathcal{I}$ into $N_\mathrm{cpu}$ sub-intervals $\mathcal{I}_n = [T_{n-1},T_n),\;n=1,\ldots,N_\mathrm{cpu}$, with $T_0=t_0$ and
$T_{N_\mathrm{cpu}}= t_{\mathrm{end}}$, according to the number  $N_\mathrm{cpu}$ of available CPUs. For an initial value problem (IVP)
\begin{align*}
	\mathbf{M}\frac{\mathrm{d}}{\mathrm{d} t} \mathbf{x}  = \mathbf{f}(\mathbf{x},t), && \mathbf{x}(t_0) = \mathbf{x}_0\;,
\end{align*}
with $\mathbf{x}:\mathcal{I}\rightarrow\mathbb{R}^{n_{\mathrm{x}}}$ and $n_{\mathrm{x}}$ the number of degrees of freedom, each
Parareal iteration $k$ solves $N_\mathrm{cpu}$ IVPs 
\begin{align}\label{eq:ivpPR}
	\mathbf{M}\frac{\mathrm{d}}{\mathrm{d} t} \mathbf{x}_n  = \mathbf{f}(\mathbf{x}_n,t), && \mathbf{x}_n(T_{n-1}) = \mathbf{X}_{n-1}^k, && t \in \mathcal{I}_n\;,
\end{align}
with $\mathbf{X}_0^k = \mathbf{x}_0$, in parallel. As the initial conditions $\mathbf{X}_{n-1}^k$ are a priori unknown
and thus do not necessarily correspond to the correct values of the continuous solution of the sequential problem,
jumps arise at the interfaces $T_{n-1}$ between the the solution of the IVP on $\mathcal{I}_{n-1}$  and the initial condition of the IVP on
$\mathcal{I}_n$.  To reduce this mismatch and converge to the original, continuous solution, 
an update is performed on each iteration \cite{Lions_2001aa}, such that for all $k$ and $n=1,\ldots,N_\mathrm{cpu}-1$
\begin{align}
	\mathbf{X}_n^k ={}& \mathcal{F}(T_n,T_{n-1},\mathbf{X}_{n-1}^{k-1})\nonumber \\
					&+ \mathcal{G}(T_n,T_{n-1},\mathbf{X}_{n-1}^{k})
					- \mathcal{G}(T_n,T_{n-1},\mathbf{X}_{n-1}^{k-1})\;. \label{eq:PRupdate}
\end{align}
Here, $\mathcal{F}(T_n,T_{n-1},\star)$ and $\mathcal{G}(T_n,T_{n-1},\star)$ are the fine and coarse operators 
that return solutions of the initial value problems \eqref{eq:ivpPR} at time $T_n$ with initial condition $\star$ at $T_{n-1}$. 
The fine propagator $\mathcal{F}$ can be computed in parallel on all subintervals $\mathcal{I}_n$ for the update formula \eqref{eq:PRupdate}. 
Therefore, it may be computationally costly to compute and thus it is 
chosen to return an accurate solution of \eqref{eq:ivpPR} by e.g. selecting a sufficiently small time step size. 
 The coarse propagator $\mathcal{G}$, however, must be executed sequentially due to the update formula \eqref{eq:PRupdate}.
 Thus it is chosen to be 
cheaper to compute (e.g. a time integrator with larger time steps) and as a consequence it is less accurate.

\subsection{Implicit/Explicit Parareal}\label{sec:impexppr}
We consider the eddy current problem \eqref{eq:eddy} and follow the approach of \cite{Dutine_2017aa} to apply an explicit time integration
scheme to the reduced ODE system that is obtained after applying the Schur complement (see \eqref{eq:ode}). To speed up the simulation
time, the Parareal algorithm is applied to the given setting. 
As it has been mentioned previously, to ensure explicit time integration schemes are stable, the time step size has to be chosen within 
the stability region of the method for the specific problem. This does not significantly affect the performance of the fine solver as, for accuracy reasons, 
smaller time steps are chosen in any case. However, for the coarse solution, a large time step size is required to ensure a reduction of simulation time
due to its sequential calling. This, however, is not possible with an explicit time integrator, as the solver would not only be inaccurate but
also unstable. 
Therefore we propose the following new approach:
\begin{itemize}
	\item On the fine level, an explicit time integration method is used with a very small time step size that ensures both accuracy as well as stability.
	\item The coarse solver performs an implicit time integration scheme with a large time step size that is less accurate but stable. 
\end{itemize}

\section{Numerical Example}\label{sec:numerics}
\begin{figure}
	\centering
	\includegraphics{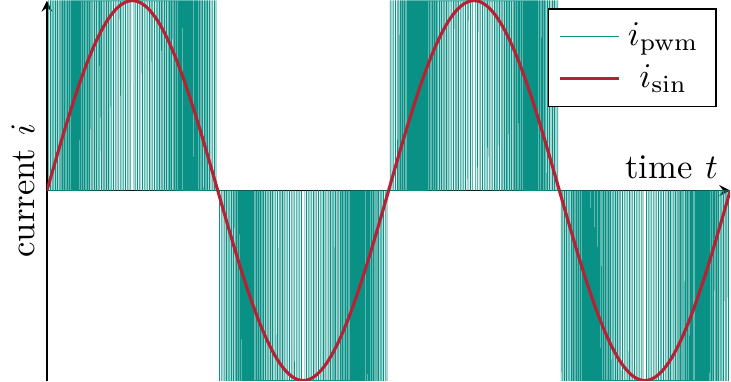}
	\caption{PWM signal $i_{\mathrm{pwm}}$ and lowest frequency component $i_{\mathrm{sin}}$.}\label{fig:pwm}
\end{figure}

To test the proposed algorithm, we solve the eddy current equation \eqref{eq:eddy} on the two dimensional model of a
single-phase isolation transformer `Mytransformer' of Figure~\ref{fig:transformer} discretised by FEMM 
\cite{Meeker_2018aa}\footnote{\url{http://www.femm.info/wiki/MyTransformer}}. 
For the given example, we set field-independent materials
and thus obtain a linear eddy current equation \eqref{eq:eddy}.
The excitation of the problem is a PWM signal $\mathbf{i}_\mathrm{pwm}(t)$ switching at frequency $f_\mathrm{pwm}$ and we 
denote by $\mathbf{i}_\mathrm{sin}(t)$
 its lowest frequency component, i.e. a $f_\mathrm{sin}$ sine wave (see Figure~\ref{fig:pwm}).

We apply the implicit/explicit Parareal algorithm for  the time interval $\mathcal{I}=(0, 0.04]\,$s with the PWM signal's
frequencies $f_\mathrm{sin}=50$\,Hz and $f_\mathrm{pwm}=10$\,kHz. 
The discretised magnetic vector potential is initialised at zero, that is $\mathbf{a}(t_0) = \bm{0}$.
The propagators are chosen as follows.

\bigskip
\textbf{Fine solution:}\\
Here the original problem is solved with the PWM excitation signal $\mathbf{i}=\mathbf{i}_\mathrm{pwm}(t)$ (see sketch in Figure~\ref{fig:pwm}).
For the time integration, the explicit Euler scheme is used as in \cite{Schops_2012aa,Dutine_2017aa} on the reduced ODE \eqref{eq:ode} 
 of the eddy current differential algebraic equation. The time step size of the method is set to $h =10^{-8}$s. 

\bigskip
\textbf{Coarse solution:}\\
On the coarse level the approach of \cite{Gander_2019aa} is followed and the eddy current problem \eqref{eq:eddy} is excited only with the lowest
frequency component of the PWM excitation $\mathbf{i}\equiv\mathbf{i}_\mathrm{sin}$. This allows using a larger time step size for the integration,
which  significantly reduces the computational cost. As described in Section~\ref{sec:impexppr}, the implicit Euler scheme is 
unconditionally stable  and a time step size of $H=0.04/N_{\mathrm{cpu}}$s is used, which corresponds to one time step per CPU.
 
\bigskip
The Parareal algorithm is iterated until the jumps of the magnetic vector potential at the interfaces between windows are  below an $l^2$ error 
with relative tolerance $\mathrm{reltol}=10^{-4}$ and absolute tolerance  $\mathrm{abstol}=10^{-10}$.

\subsection{Simulation Results}
\begin{figure}[t]
    \centering
    \includegraphics{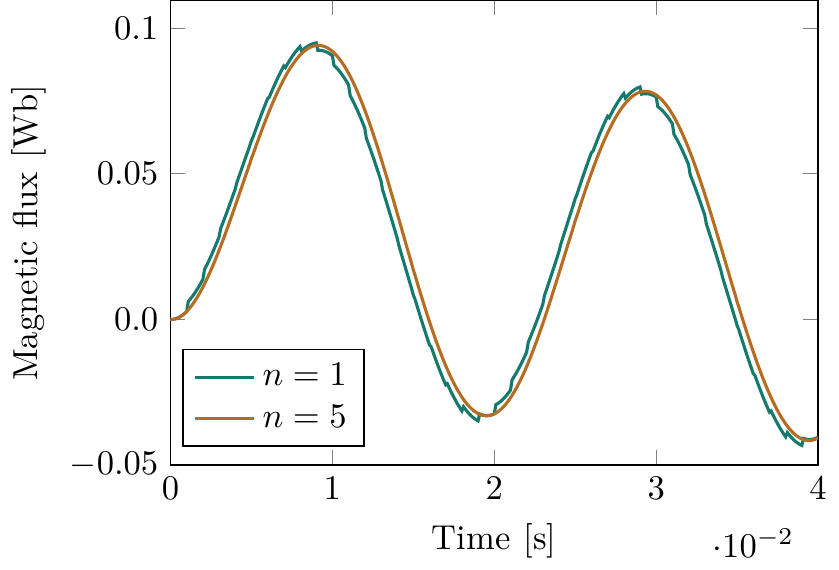}
    \caption{Results of Parareal after $n=1,5$ iterations.}
    \label{fig:results}
\end{figure}

For a number of $N_\mathrm{cpu}=40$ CPUs, the Parareal algorithm converges to the required tolerance after $n=5$ iterations.
This implies a theoretical speed-up of $40/5=8$ with respect to a sequential explicit computation,
when neglecting communication costs and coarse level computations. 
Using Intel Dualcore i7 (3,1 GHz) based hardware and an implementation in Matlab 2020b, the evaluation of the 
coarse propagator (implicit Euler) on the overall time domain requires less than $1$s, while each fine propagator call 
(explicit Euler, PWM right-hand-side) runs approximately $30$s. 
The resulting magnetic flux for the first iteration ($n=1$) and after the algorithm is converged is
given in Figure~\ref{fig:results}. Note the artificial jumps in the solution for $n=1$ that are smoothed out in the 5th iteration's solution.

\section{Summary and Outlook}\label{sec:summary}
This article proposes a combination of explicit and implicit time integration methods for the parallel-in-time method Parareal. Its application
is exemplified with the eddy current equation solved for a simple model of a single-phase isolation transformer. 
The result confirms the fast convergence of the
algorithm, which yields in the best case for 40 processors a theoretical speed up of  8 in comparison to explicit methods, which has been shown 
to be faster than implicit methods in comparable studies \cite{Dutine_2017aa}.

In future work nonlinearities should be included in the model to exploit all the advantages of the usage of explicit time integration methods
for the eddy current equation (see \cite{Kaehne_2020aa}). Furthermore, more evolved explicit time integration schemes such as the Runge-Kutta-Chebyshev
method are to be studied \cite{Kaehne_2020aa}. Hereby the stability region of the method is increased and larger time step sizes can be used.

\section*{\small Acknowledgement}
{\footnotesize This work is supported by the Graduate School CE within the Centre for Computational 
	Engineering at Technische Universit\"at Darmstadt and DFG Grants 
	SCHO1562/1-2, CL143/11-2 and BMBF Grant 05M2018RDA (PASIROM). \par} 
\renewcommand*{\bibfont}{\small}
\printbibliography

\end{document}